\documentclass{jfl}
\usepackage{enumerate}
\usepackage{enumerate}
\usepackage{graphicx}
	
\newtheorem*{theoetoile}{Theorem} 
\newtheorem*{coroetoile}{Corollary} 

\newcommand{\rsub}[2]{#1^{[#2]}}

\newcommand{\mC}{\mathcal{C}}
\newcommand{\mD}{\mathcal{D}}

\newcommand{\mI}{\mathcal{I}}

\newcommand{\mL}{\mathcal{L}}

\newcommand{\mP}{\mathcal{P}}

\newcommand{\pf}{\textit{Proof: }}
\newcommand{\dom}{\text{dom}}

\title{Locally constant functions in $C$-minimal structures}
\author{Pablo Cubides Kovacsics}
\thanks{The author was supported by the Marie Curie Initial Training Network in Mathematical Logic - MALOA
- From MAthematical LOgic to Applications, PITN-GA-2009-238381. Many thanks to the referee for his/her careful reading of the article.}

\begin{document}

\parindent=0in

\maketitle

\begin{abstract}
Let $M$ be a $C$-minimal structure and $T$ its canonical tree (which corresponds in an ultrametric space to the set of closed balls with radius different than $\infty$ ordered by inclusion). We present a description of definable locally constant functions $f:M\rightarrow T$ in $C$-minimal structures having a canonical tree with infinitely many branches at each node and densely ordered branches. This provides both a description of definable subsets of $T$ in one variable and analogues of known results in algebraically closed valued fields.
\end{abstract}

\

\

\

This paper studies the behavior of definable locally constant functions $f:M\rightarrow T$ where $M$ is a $C$-minimal structure, $T$ is its canonical  tree, each node of $T$ has infinitely many branches and branches are densely ordered. These conditions are fulfilled in different $C$-minimal ultrametric spaces such as algebraically closed valued fields, where $T$ corresponds to the set of closed balls not having $\infty$ as their radius (i.e., not being a singleton) ordered by inclusion. Locally constant functions appear naturally when one gives a general description of definable sets in $C$-minimal structures such as the cell decomposition theorem proved by Haskell and Macpherson in \cite{macphersonETAL:94}, especially if the theory satisfies the exchange property. In particular, we show that (all terms to be later defined):

\begin{theoetoile} Let $f:M\rightarrow T$ be a partial definable locally constant function. Then $\dom(f)$ can be decomposed (possibly adding parameters) into cells $D_1,\ldots,D_n$ such that for each $1\leq i\leq n$ one and only one of the following conditions holds:
	\begin{enumerate}
	\item $f(D_i)$ is an antichain.
	\item $f(D_i)$ is a chain.
\end{enumerate}
\end{theoetoile}

We use the result to show that the domain of a partial definable function from $M$ to a branch $B$ of $T$ can be decomposed into finitely many cells $D_1,\ldots,D_n$ where the function can be factored as in the following diagram (modulo some finite set in the range)

\begin{center}
	\includegraphics[width=35mm]{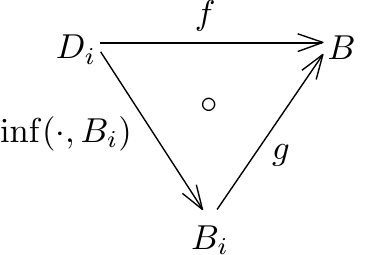}
\end{center}

where $B_i$ is a branch of $T$ and the function $\inf(\cdot,B_i)$ sends $\alpha$ to the infimum in $T$ of $\alpha$ and $B_i$. An analogue of this result is given for functions having their range in the set of cones at a given point in $T$. Finally, the theorem serves also to give a description of definable subsets of $T$ in one variable and to prove the following two results on $C$-minimal expansions of algebraically closed valued fields: 

\begin{coroetoile}
Let $M$ be a definably complete $C$-minimal expansion of an algebraically closed valued field $K$. We denote by $vK$ the valuation group. Let $f:K\rightarrow vK$ be a partial definable function. Then there are a finite set $F\subseteq f(K)$ and $\beta_1,\ldots,\beta_n$ such that $\dom(f)\setminus f^{-1}(F)$ can be decomposed into cells $D_1\ldots,D_n$ satisfying that for all $1\leq i\leq n$ and all $a\in f(D_i)$ there is $a_i\in vK$ such that
	$$\{x\in D_i: f(x)=a\}=\{x\in D_i: v(x-\beta_i)=a_i\}$$
\end{coroetoile}

\begin{coroetoile}
Let $M$ be a $C$-minimal expansion of an algebraically closed valued field $K$. We denote by $K/v$ the residue field. Let $f:K\rightarrow K/v$ be a partial definable function. Then, there are a finite set $F\subseteq f(K)$, a positive integer $s$, elements $\alpha_1,\ldots,\alpha_s,\beta_1,\ldots,\beta_s\in K$ and definable functions $h_{i}:K/v\rightarrow K/v$ for $1\leq i\leq s$ such that for $X:=\dom(f)\setminus f^{-1}(F)$ and $t_i:=v(\alpha_i-\beta_i)$ (for $1\leq i\leq s$)
 $$X\subseteq \bigcup_{1\leq i\leq s} \Lambda_{t_i}(\alpha_i),$$
(where $\Lambda_{t_i}(\alpha_i)$ is the closed ball of radius $t_i$ centered at $\alpha_i$, see Section \ref{sec1} for more details) and for $1\leq i\leq s$ and $x\in \Lambda_{t_i}(\alpha_i)\cap X$,
$$f(x)=h_{i}\left(\left(\frac{x-\alpha_i}{\alpha_i-\beta_i}\right)/v\right).$$
\end{coroetoile}

\

Section 1 contains a brief introduction to $C$-minimality together with the definitions and properties later needed. In Section 2 we prove the above mentioned theorem and its corollaries. Finally, in Section 3  the description of definable subsets of $T$ is presented. For a more detailed introduction to $C$-minimality the reader is invited to see any of \cite{macphersonETAL:94, macphersonETAL:96, delon:11, cubides:12}.

\section{$C$-minimality}\label{sec1}

\begin{definition} Let $C(x, y, z)$ be a ternary relation. A \emph{$C$-set} is a structure $(M,C)$
satisfying axioms (C1)-(C4):

\begin{enumerate}
\item[(C1)] $\forall xyz (C(x, y, z)\rightarrow C(x, z, y))$
\item[(C2)] $\forall xyz (C(x, y, z)\rightarrow \neg C(y, x, z))$
\item[(C3)] $\forall xyzw (C(x, y, z)\rightarrow(C(w, y, z)\vee C(x, w, z)))$
\item[(C4)] $\forall xy (x \neq y \rightarrow C(x, y, y))$
\item[(D)] 	$\exists x,y(x\neq y) \wedge \forall xy (x \neq y \rightarrow \exists z(z\neq y \wedge C(x,y,z)))$
\end{enumerate}	

If in addition $(M,C)$ satisfies axiom (D) we say it is a \emph{dense} $C$-set.
\end{definition}

Throughout the paper a tree $T$ is a partially ordered set such for all $a\in T$ the set $a_<:=\{b\in T: b<a\}$ is totally ordered (not necessarily well-ordered) and any two elements have a common lower bound. Every $C$-set $M$ interprets a tree denoted $T(M)$ which is a meet semi-lattice (every two nodes $a,b\in T(M)$ have an infimum denoted $\inf(a,b)$) such that for every node there is a leaf above it. In fact the class of such trees (called \emph{good trees}) is bi-interpretable with the class of $C$-sets (see \cite{delon:11}). Following Delon in \cite{delon:11}, we identify the set of leaves in $T(M)$ with $M$. We denote by $T$ the set of elements in $T(M)$ which are not a leaf. A $C$-structure is simply a $C$-set with possibly extra structure. In what follows we work in a fixed $C$-structure $M$ in a language $L$. We use lower case Greek letters $\alpha,\beta,\gamma$ to denote both elements of $M$ and leaves in $T(M)$ and lower case letters $a,b,c$ to denote arbitrary elements in $T(M)$. For $a\in T$, we define an equivalence relation $E_a$ on $a_{>}$ (i.e. $\{b\in T(M): b>a\}$) by

$$x E_a y \Leftrightarrow a< \inf(x,y).$$

Equivalence classes are called \emph{cones at $a$}. \emph{The branching number of $a$}, denoted by $bn(a)$, is the number of equivalence $E_a$-classes. For $a,b\in T$ such that $a<b$, \emph{the cone of $b$ at $a$}, denoted by $\Gamma_a(b)$, is the $E_a$-class of $b$. We abuse notation using $\Gamma_a(b)$ to denote also the subset of $M$ defined by $\Gamma_a(b)\cap F$, where $F$ is the set of leaves in $T(M)$ (again identifying this set with $M$). In particular, for $\alpha\in M$ and $a\in T$ such that $a<\alpha$, the cone $\Gamma_a(\alpha)$ will be usually taken to be the set $\{\beta\in M: a<\inf(\alpha,\beta)\}$. For $\alpha,\beta\in M$, we then have that 
$$\Gamma_{\inf(\alpha,\beta)}(\alpha)=\{x\in M:M\models C(\beta,x,\alpha)\}.$$
For practical reasons we treat $M$ as a cone at $-\infty$, that is, we extend $T(M)$ by adding a new element $-\infty$ satisfying $-\infty<a$ for all $a\in T(M)$ and we let $\Gamma_{-\infty}(\alpha):=M$ for all $\alpha\in M$. For $a\in T(M)$, an \emph{$n$-level set at $a$} is the set $\{x\in T: a\leq x\}$ with $n$ cones at $a$ removed, provided that if $n\geq 1$, then $bn(a)>n$. For $b_1,\ldots,b_n\in T$ such that $a<b_i$ and $\neg b_iE_a b_j$ for all $1\leq i<j\leq n$, the expression $\Lambda_a(b_1,\ldots,b_n)$ denotes the $n$-level set where the $n$ cones removed are $\Gamma_a(b_i)$ for $1\leq i\leq n$. In particular, for $a\in T(M)$, the $0$-level set at $a$ is denoted in symbols by $\Lambda_a$. If $a\in T$, $\Lambda_a$ corresponds to the union of all cones at $a$; if $a=\alpha\in M$, then $\Lambda_\alpha=\{\alpha\}$. For $\alpha,\beta\in M$ we have that
$$\Lambda_{\inf(\alpha,\beta)}=\{x\in M:M\models\neg C(x,\alpha,\beta)\}.$$
Both for cones and $n$-level sets, the point $a$ is called its \emph{basis} and we let $\min(\cdot)$ to be the function sending cones and $n$-level sets to their bases. Finally, for $a,b\in T(M)\cup\{-\infty\}$, such that $a<b$, \emph{the interval} $(a,b)$ denotes in $T(M)$ the set $\{x\in T: a<x<b\}$ and in $M$ the set $\Gamma_a(b)\setminus \Lambda_b$. Cones, intervals and $n$-level sets can be seen as subsets of $T(M)$ or $M$ and we usually let the context specify which one is intended. We denote the set of all cones (including $M$) by $\mC$, the set of all intervals by $\mI$ and the set of $n$-levels by $\mL_n$ for $n<\omega$. For $a\in T$, we let $\mC[a]$ be the set of cones at $a$. In $M$, the set of cones $\mC$ forms a uniformly definable basis of clopen sets (``uniformly'' here means the same formula is used, changing its parameters, to define all basic open sets ). We work with the topology generated by this basis which is Hausdorff and totally disconnected.

\begin{definition} A $C$-structure $M$ is \emph{$C$-minimal} if for every elementary equivalent structure $N\equiv M$, every definable subset $D\subseteq N$ is definable by a quantifier free formula using only the predicate $C$. A complete theory is $C$-minimal if it has a $C$-minimal model.
\end{definition}

From now on $M$ will be a dense $C$-minimal structure. For $\alpha\in M$, the branch of $\alpha$ is the set $Br(\alpha):=\{a\in T(M):a\leq \alpha\}$. $C$-minimality implies that $Br(\alpha)$ is \emph{o-minimal} in the sense that any subset of $Br(\alpha)$ definable in $M$ is a finite union of intervals and points. Analogously, again by $C$-minimality, for each $a\in T(M)$, the set of all cones at $a$ is strongly minimal in the sense that for any elementary extension $N$ of $M$, any subset of cones at $a$ which is definable in $N$ is either finite or cofinite (see \cite{macphersonETAL:94, cubides:12}). Compare the previous statement to the fact that in an algebraically closed valued field, the value group is o-minimal and the residue field is algebraically closed, hence strongly minimal (see \cite{haskellETAL:08}). In analogy with o-minimal theories,  Haskell and Macpherson proved in \cite{macphersonETAL:94} a cell decomposition theorem for dense $C$-minimal structures. The definition of a cell is more involved than in o-minimality and for the purpose of this article, the reader may think of 1-cells as finite unions of either points, cones, intervals or $n$-level sets (formal conditions defining a 1-cell can be found in \cite{macphersonETAL:94} and \cite{cubides:12}). For a set $A$, we use $\rsub{A}{r}$ to denote the set of all subsets of $A$ of cardinality $r$. Thus 1-cells are subsets \emph{of type} $\rsub{Z}{r}$ for $r$ a positive integer and $Z$ either $M,\mC,\mI$ or $\mL_n$ for $n<\omega$. The following proposition corresponds to the cell decomposition theorem for formulas in one variable (this is Lemma 2.9 in \cite{macphersonETAL:94}; a proof of the proposition as here stated can be found in \cite{cubides:12}, Proposition 22.)

\begin{proposition}\label{celldecomposition1}
Let $\phi(x,y)$ be an $L$-formula with $|x|=1$. There is a finite definable partition $\mP$ of $M^{|y|}$ such that for each $A\in \mP$ there are $L$-formulas $\psi_{1}^A(x,y),\ldots,\psi_{n_A}^A(x,y)$ satisfying that whenever $\alpha\in A$, $\{\psi_{1}^A(M,\alpha),\ldots,\psi_{n_A}^A(M,\alpha)\}$ is a 1-cell decomposition of $\phi(M,\alpha)$. Moreover, each formula $\psi_{j}^A(x,y)$ defines the same type of 1-cell for each $\alpha\in A$.
\end{proposition}

We state a lemma from \cite{macphersonETAL:94} (fact 1) which will be later used (a proof can be found in the appendix of \cite{cubides:12}, Lemma 9).

\begin{lemma}\label{fait1} Let $S\subseteq M$ be a definable set. Then, there is no $\alpha\in M$ such that for an infinite number of nodes $a\in Br(\alpha)$ we have both
\begin{equation*}
\label{eq:1}\tag{$\ast$} \Lambda_a(\alpha)\cap S\neq \emptyset \text{ and } \Lambda_a(\alpha)\cap (M\setminus S)\neq \emptyset .
\end{equation*}
\end{lemma}

We turn now our attention to definable functions. The study of definable functions in dense $C$-minimal structures started in \cite{macphersonETAL:94} and was mainly used to prove the cell decomposition theorem. For a definable function $f:M\rightarrow T$ their study shows that $\dom(f)$ can be definably partitioned into three sets $F\cup I\cup K$. where $F$ is finite, $f$ is locally constant on $K$ and on $I$, locally, the range of $f$ behaves like a $C$-structure and $f$ behaves like an isomorphism of $C$-sets. This can be made precise by defining first a $C$-relation on any subset $A\subseteq T$ which is an antichain as follows: let $T'$ denote the reduct of $T$ to the language $L_{\inf}:=\{\leq,\inf\}$ of meet semi-lattices and let $T[A]:=\left\langle A\right\rangle_{T'}$ be the $L_{\inf}$-substructure generated by $A$ (notice that the universe of $T[S]$ is the $\inf$-closure (denoted $cl_{\inf}(A)$) of $A$ and that $A$ itself corresponds to the set of leaves in $T[A]$); it is easy to show that $A$, which from now on we denote by $M[A]$, is a $C$-set with canonical tree $T[A]$.

\begin{definition} Let $M$ and $N$ be two $C$-sets. A function $f:M\rightarrow N$ is called a $C$-isomorphism if it is injective and preserves the $C$-relation and its negation. For $T$ the canonical tree of $M$ and $A\subseteq T$ an antichain, a function $f:M\rightarrow A$ is a $C$-isomorphism if $f:M\rightarrow M[A]$ is a $C$-isomorphism.
\end{definition}

The previous result of Haskell and Macpherson in \cite{macphersonETAL:94} corresponds then to (Theorem 3.11 in \cite{macphersonETAL:94}):

\begin{proposition}[Haskell-Macpherson]\label{monotonicity}
	Let $f:M\rightarrow T$ be a partial definable function. Then $\dom(f)$ can be decomposed into finitely many 1-cells such that on each infinite cell, $f$ is either locally constant or a local $C$-isomorphism.
\end{proposition}

Also in \cite{macphersonETAL:94}, Haskell and Macpherson characterized those $C$-minimal structures for which the algebraic closure lacks the exchange property by the existence of ``bad functions'', which are definable functions $f:M\rightarrow T$ for which there is a cone $D$ such that $f\upharpoonright D$ is a $C$-isomorphism. 

\begin{theorem}[Haskell-Macpherson]\label{badfunction} $Th(M)$ has the exchange property if and only if $M$ has no bad function.
\end{theorem}

Thus, in a $C$-minimal structure $M$ such that $Th(M)$ satisfies the exchange property, all definable functions $f:M\rightarrow T$ are locally constant.

\section{Locally constant functions to the canonical tree}

Through this section we work in a dense $C$-minimal structure $M$ where for all $a\in T$, $bn(a)\geq \aleph_0$ and each branch $B\subseteq T$ is densely ordered. For $a,b\in T$ we use the relation symbol $a \parallel b$ to say that $a$ and $b$ are incomparable, and $a\nparallel b$ to say they are comparable. We start with some preparation lemmas. 

\begin{lemma}\label{denseanti}
	Let $A\subseteq T$ be a definable antichain. Then the set 
	$$A_0:=\{a\in A: a\text{ has a predecessor in } cl_{\inf}(A)\},$$ is finite (by `predecessor' we mean `immediate predecessor').  
\end{lemma}

\pf Given that $A$ is definable so is $A_0$. Suppose towards a contradiction that $A_0$ is infinite. For each $a\in A_0$ let $a^-$ be its predecessor in $cl_{\inf}(A)$. Then consider the definable set $D:=\bigcup_{a\in A_0} \Lambda_a$ and the tree 
$$T(D):=\{a\in T: \exists\beta(\beta\in D\wedge a<\beta) \wedge \exists\beta(\beta\in M\setminus D\wedge a<\beta)\}.$$
By Proposition 3.7 in \cite{delon:11}, $T(D)$ has finitely many branches (\emph{i.e.} finitely many nodes $a$ such that $bn(a)\geq 2$, where the function $bn$ is calculated in $T(D)$). We show that the set $A':=\{a^-:a\in A_0\}$ is included in $T(D)$ which gives a contradiction since $A'$ is an infinite antichain. Clearly, for each $a\in A'$ there is $\beta\in D$ such that $a<\beta$. By the density assumption on the branches, there is $\beta\in \gamma_{a^-}(a)\setminus \Lambda_a$, hence $\beta\notin D$. $\square$  

\

The previous lemma shows in particular that if $A$ is an infinite antichain, there is an sub-antichain $A'\subseteq A$ such that $|A\setminus A'|$ is finite and $M[A']$ is a dense $C$-minimal structure. 

\begin{lemma}\label{antichain1}
Let $a\in T$ and let $f:\Lambda_a\rightarrow T$ be a partial definable function such that $f\upharpoonright\Gamma_a(\alpha)$ is constant for every $\alpha\in \Lambda_a$. Then there is a finite subset $F\subseteq f(\Lambda_a)$ such that $f(\Lambda_a)\setminus F$ is an antichain.
\end{lemma}

\pf Given that for every $\alpha\in \dom(f)$, $f\upharpoonright\Gamma_a(\alpha)$ is constant, we can treat $f$ as a partial function from $\mC[a]$ to $T$. If $f(\Lambda_a)$ is finite there is nothing to prove, so we may assume that $f$ has infinite range. By compactness and $C$-minimality, one of the following happens (possibly in an elementary extension of $M$): either $f(\Lambda_a)$ contains an infinite definable chain or there is a positive integer $k$ such that every chain contained in $f(\Lambda_a)$ has cardinality at most $k$. The former is impossible since adding parameters to define this chain we get a definable function from a strongly minimal set to an infinite ordered set. Thus suppose there is a positive integer $k$ as above mentioned. This implies already the result since if there were no finite subset $F\subseteq f(\Lambda)$ such that $f(\Lambda)\setminus F$ is an antichain, the set of cones at $a$ for which their image under $f$ is maximal in $f(\Lambda)$ would be an infinite coinfinite set, which contradicts strong minimality. $\square$

\begin{lemma}\label{antichain2}
Let $A\subseteq T$ be an infinite antichain and $\{\Gamma_i:i<\kappa\}$ be an infinite set of cones satisfying
\begin{enumerate}
	\item for all $i<\kappa$ the basis of $\Gamma_i$ is an element of $A$.
	\item there is a positive integer $N$ such that for every $a\in A$ there are at most $N$ cones in $\{\Gamma_i:i<\kappa\}$ for which $a$ is the basis.
\end{enumerate}
Then $H=\bigcup \{\Gamma_i:i<\kappa\}$ is not a definable set.
\end{lemma}

\pf Suppose $H$ is definable. By $C$-minimality there is a decomposition into 1-cells $\mD=\{D_1,\ldots,D_s\}$ of $H$.

\

\textbf{Claim 1:} For each definable set $D\subseteq M$ there are finitely many cones in $\{\Gamma_i:i<\kappa\}$ intersecting both $D$ and $M\setminus D$.

\

\textit{Proof of Claim 1:} Suppose towards a contradiction this is not the case. Consider the tree
$$T(D):=\{a\in T: \exists\beta(\beta\in D\wedge a<\beta) \wedge \exists\beta(\beta\in M\setminus D\wedge a<\beta)\}.$$
Since $A$ is an antichain, there are infinitely many points $a\in T(D)$ such that $bn(a)\geq 2$ (one for each cone in $\{\Gamma_i:i<\kappa\}$ intersecting both $D$ and $M\setminus D$). But again, by Proposition 3.7 in \cite{delon:11}, $T(D)$ must have finitely many branches, a contradiction.  This completes the claim.

\

By the claim there must be a 1-cell $D\in \mD$ that contains infinitely many cones in $\{\Gamma_i:i<\kappa\}$. If $D$ is a 1-cell of type $\rsub{Z}{r}$ where $Z$ is any of $\mC, \mI$ or $\mL_n$ for $n<\omega$ and $r>1$, by the claim there is a 1-cell $E$ of type $\rsub{Z}{1}$ contained in $D$ and containing infinitely many cones in $\{\Gamma_i:i<\kappa\}$. If $E$ is a cone or an $n$-level set, let $b$ be its basis. If $E$ is an interval, let $b$ be its left-ending point. Condition (2) together with the assumption that there are infinitely many cones in $\{\Gamma_i:i<\kappa\}$ contained in $E$ implies that there is $i_0<\kappa$ such that $\Gamma_{i_0}$ is contained in $E$ and $b<a:=\min(\Gamma_{i_0})$. Since $bn(a)\geq\aleph_0$, condition (2) and the claim entail that there is $\beta\in E$ such that $a<\beta$ and the unique cone in $\{\Gamma_i:i<\kappa\}$ containing $\beta$, say $\Gamma_{i_1}$, is contained in $E$ and has a different basis than $a$. But then both $\min(\Gamma_{i_0})$ and $\min(\Gamma_{i_1})$ belong to $Br(\beta)$ which contradicts that $A$ is an antichain by condition (1). $\square$

\begin{theorem}\label{localcst}
	Let $f:M\rightarrow T$ be a partial definable and locally constant function. Then $\dom(f)$ can be decomposed (possibly adding parameters) into 1-cells $D_1,\ldots,D_n$ such that for each $1\leq i\leq n$ one of the following conditions holds:
	\begin{enumerate}
	\item $f(D_i)$ is an antichain.
	\item $f(D_i)$ is a chain.
\end{enumerate}
\end{theorem}

\pf It is enough to show that there is finite definable partition $\mP$ of $\dom(f)$ such that each $D\in \mP$ has the following property  

\begin{equation}\label{propertyP}
  \tag{\textbf{P}} f(D) \text{ is either a finite union of antichains or a finite union of chains.}
\end{equation}

For suppose $\mP$ is such a partition. Fix $D\in \mP$ and suppose furthermore that $f(D)$ is a finite union of antichains. Define subsets $(X_n)_{n<\omega}$ of $f(D)$ by 

\begin{enumerate}
	\item[$\bullet$] $X_0$ is the set of maximal elements in $f(D)$;
	\item[$\bullet$] $X_{i+1}$ is the set of maximal elements in $\displaystyle f(D)\setminus \bigcup_{j=1}^i X_j$.
\end{enumerate}

By construction, for all $n<\omega$ each $X_n$ is either empty or a definable antichain. Moreover for $i<j<\omega$, $X_i\cap X_j=\emptyset$. Our assumption implies there is $n<\omega$ such that $X_{m}=\emptyset$ for all $m>n$. Taking $f^{-1}(X_j)\cap D$ we get a partition of $D$ into finitely many pieces each corresponding to the preimage of an antichain. Now suppose that $f(D)$ is a finite union of chains. Possibly adding new parameters, we can partition this finite union of chains in a definable way getting finitely many disjoint definable chains, say $I_1,\ldots I_n$. One can then partition $D$ into subsets $f^{-1}(I_i)\cap D$ for each chain $1\leq i\leq n$. The image of each subset $f^{-1}(I_i)\cap D$ is clearly a chain. The result is obtained by decomposing into cells all these partitions for each $D\in\mP$. 

\

In addition, given $\mP$ a finite definable partition of $\dom(f)$, it suffices to show property \ref{propertyP} for each $D\in \mP$ modulo a finite subset $F_D\subseteq D$. This is so because in $\mP':=\{D\setminus F_D: D\in \mP\}\cup \{F_D:D\in \mP\}$, each finite subset $F_D$ already satisfies property \ref{propertyP} ($f(F_D)$ is is in particular a finite union of antichains). 

\

Let $g:\dom(f)\rightarrow T$ be the definable function sending $\alpha$ to the basis of the maximal cone $\Gamma\subseteq\dom(f)$ containing $\alpha$ on which $f$ is constant. By $C$-minimality, $g$ is well-defined given that $f$ is locally constant. Consider the following formula:

$$\phi(x,y):= \Gamma_{g(y)}(x)\subseteq\dom(f) \wedge \forall z \in \Gamma_{g(y)}(x) (f(z)=f(y)).$$

For $\alpha\in\dom(f)$, the set $\phi(M,\alpha)$ is simply $\Lambda_{g(\alpha)}\cap f^{-1}(f(\alpha))$, the union of all cones at $g(\alpha)$ contained in $f^{-1}(f(\alpha))$. By definition, $\phi(M,\alpha)$ is a union of cones at $g(\alpha)$. By $C$-minimality (in particular by \ref{celldecomposition1}) there are positive integers $s$ and $n$ such that $\phi(M,\alpha)$ is either the union of at most $s$ cones at $g(\alpha)$ or an $m$-level at $g(\alpha)$ for some $m\leq n$. As a consequence, there is a finite definable partition $\mP$ of $\dom(f)$ such that for each $A\in \mP$, all $\alpha\in A$ satisfy that $\phi(M,\alpha)$ is of the same type, i.e., either a union of exactly $r$ cones or an $m$-level set with $r\leq s$ and $m\leq n$. Notice that $\phi(\beta,\alpha)$ implies that $\phi(M,\beta)=\phi(M,\alpha)$, which shows that if $\alpha\in A$ then $\phi(M,\alpha)\subseteq A$. By the initial remark, it suffices to prove that each subset $A\in \mP$ can be definably partitioned into finitely many sets satisfying property \ref{propertyP}. Fixing $A\in \mP$, we split in cases depending on the form of $\phi(M,\alpha)$:

\

\textbf{Case 1:} Suppose that $\phi(M,\alpha)$ is of type $\rsub{\mC}{r}$. Let us assume towards a contradiction that $f(A)$ contains an infinite chain $I$. In particular without loss of generality we may assume that $I$ is definable. Consider the definable set $D:=f^{-1}(I)\cap A$ (possibly adding new parameters for $I$ if necessary) and, for the sake of the argument, take $f=f\upharpoonright D$. For $B:=\{g(\alpha)\in T: \alpha\in D\}$ we have either that

\begin{enumerate}
	\item[(I)] $B$ is a finite union of antichains,
	\item[(II)] $B$ contains an infinite definable chain (possibly replacing $M$ by an elementary extension).
\end{enumerate}

Suppose that (I) is the case. Let $B_0$ be the antichain of maximal elements in $B$. By compactness, either there is $a\in B_0$ such that infinitely many cones at $a$ have different images in $I$ or there is an integer $s$ such that for each a $a\in B_0$ there are less than $s$ cones at $a$ having different images in $I$. If the former we get a contradiction with Lemma \ref{antichain1} (notice that since $a\in B_0$, the function $f\upharpoonright (D\cap \Lambda_a)$ is constant on cones at $a$). If the latter, since $bn(a)\geq\aleph_0$ for each $a\in T$, we get a contradiction with Lemma \ref{antichain2}. So suppose (II) is the case and let $\delta$ (possibly from an elementary extension) be such that there is $J$ an infinite definable chain contained in $B\cap Br(\delta)$. Again by Lemma \ref{antichain1}, there is no node $a\in J$ such that $f(\Lambda_a)$ is infinite. But then we have that $f^{-1}(I)\cap \Lambda_a(\delta)\neq\emptyset$ and $(M\setminus f^{-1}(I))\cap \Lambda_a(\delta)\neq\emptyset$ for all $a\in J$, which contradicts Lemma \ref{fait1}. Notice that here we use again that $bn(a)\geq\aleph_0$. This shows that $f(A)$ is a finite union of antichains (so $A$ satisfies property \ref{propertyP}).

\

\textbf{Case 2:} Assume $\phi(M,\alpha)$ is a 0-level set. By definition of $g$ and the Case 2 assumption, the set $B:=\{g(\alpha)\in T: \alpha\in A\}$ is an antichain. Suppose towards a contradiction that $f(A)$ contains an infinite definable chain $I$. Let $B_0:=\{g(\alpha):\alpha\in f^{-1}(I)\cap A\}$. By \ref{denseanti} we can suppose without loss of generality that $M[B_0]$ is dense. The function $h:M[B_0]\rightarrow I$ defined by $h(g(\alpha))=f(\alpha)$ contradicts \ref{monotonicity} given that $h$ has infinite range but can be neither a $C$-isomorphism on a cone ($I$ is totally ordered) nor a constant on a cone since we would contradict the minimality in the definition of $g$. Therefore, by compactness, $f(A)$ must be a finite union of antichains. 

\

\textbf{Case 3:} Suppose that $\phi(M,\alpha)$ is a 1-level set and let $B$ be defined as in Case 2. By compactness and Lemma \ref{antichain2}, there is a positive integer $k$ such that any definable subset $B'\subseteq B$ which is an antichain has cardinality less than $k$ (here the cones we use in order to apply Lemma \ref{antichain2} are the omitted cones on each 1-level set $\phi(M,\alpha)$). This implies that $B$ can be partitioned into finitely many totally ordered sets $B_1,\ldots,B_k$. Thus, without loss of generality we may suppose that $B$ is totally ordered. Since $bn(a)\geq\aleph_0$, the function $\hat{f}:B\rightarrow T$ defined by $\hat{f}(a)=f(\alpha)$ for some (all) $\alpha$ such that $g(\alpha)=a$ is well-defined. We study the behavior of $\hat{f}$. Consider the subsets of $B$
\small
\begin{center}
	$X:=\{b\in B: \text{there is an interval $I\subseteq B$ s.t. } (b\in I \wedge \forall a,a'\in I(\hat{f}(a)\nparallel \hat{f}(a')))\}$,
\
	$Y:=\{b\in B: \text{there is an interval $I\subseteq B$ s.t. } (b\in I \wedge \forall a,a'\in I(\hat{f}(a)\parallel \hat{f}(a')))\}$.
	\end{center}
\normalsize
We show that $B\setminus (X\cup Y)$ is a finite set. For if not, by o-minimality there is an infinite interval $J$ contained in $B\setminus (X\cup Y)$. But since $J$ is infinite, for every $N<\omega$ there is $c\in J$ such that for

\begin{enumerate}
	\item[] $X_c:=\{d\in J: \hat{f}(c)\nparallel \hat{f}(d)\}$,
	\item[] $Y_c:=\{d\in J: \hat{f}(c)\parallel \hat{f}(d)\}$,
\end{enumerate}

we have that $|X_c| > N$ is discrete which by compactness contradicts o-minimality using the density assumption on the branches of $T$. Therefore without loss of generality we may assume that $B=X\cup Y$. We show that $\hat{f}(X)$ is a finite union of chains. Suppose for a contradiction that there are infinitely many pairwise incomparable chains in $\hat{f}(X)$. The set of left endpoints of maximal intervals $I_0\subseteq X$ such that $\hat{f}(I_0)$ is totally ordered in $\hat{f}(X)$, will be a discrete infinite definable subset of $B$, which contradicts o-minimality. Therefore $\hat{f}(X)$ is a finite union of chains. Analogously, suppose that $\hat{f}(Y)$ is an infinite union of antichains. Then, the set of left end points of maximal intervals $I_0\subseteq Y$ such that $\hat{f}(I_0)$ is an antichain will be an infinite definable discrete subset of $B$, contradicting again o-minimality of a branch of $T$. Hence $\hat{f}(Y)$ is a finite union of antichains. This shows that both $X$ and $Y$ have property \ref{propertyP}. 

\

\textbf{Case 4:} Suppose that $\phi(M,\alpha)$ is an $n$-level set for $n>1$. As in case 1, let us assume towards a contradiction that $f(A)$ contains an infinite definable chain $I$, let $D:=f^{-1}(I)\cap A$ and $B:=\{g(\alpha)\in T: \alpha\in D\}$. Again we have either that:

\begin{enumerate}
	\item[(I)] $B$ is a finite union of antichains.
	\item[(II)] $B$ contains an infinite definable chain.
\end{enumerate}

If (I) is the case we get the same contradiction as in Case 1. So suppose (II) and let $J$ be an infinite definable chain contained in $B$ and $\delta$ an element (again possibly in an elementary extension) such that $Br(\delta)$ contains $J$. For $a\in J$, given that $bn(a)\geq\aleph_0$, we have $f^{-1}(I)\cap \Lambda_a(\delta)\neq\emptyset$. Moreover, since $n\geq 2$ there is at least one cone at $a$ which is not in $f^{-1}(I)\cap \Lambda_a(\delta)$, hence $M\setminus (f^{-1}(I)\cap \Lambda_a(\delta))\neq\emptyset$ contradicting Lemma \ref{fait1}. Therefore we have that $f(A)$ is a finite union of antichains. $\square$

\begin{definition}\label{incochain} Let $\mathfrak{C}$ be a finite set of definable chains of $T$. For $C\in \mathfrak{C}$ let $lp(C)$ be its left endpoint (possibly $-\infty$). Let 
$$A=\{\inf(lp(C_0),lp(C_1)): C_0\neq C_1\in \mathfrak C\}.$$ The set $\mathfrak C$ is called a set of \emph{incomparable chains} if the following conditions hold

\begin{enumerate}[(a)]
	\item for all $C_0\neq C_1\in \mathfrak C$, all $c\in C_0$ and all $d\in C_1$ we have that $c||d$;
	\item $A$ is an antichain and $Aut(M[A])$ acts transitively on $M[A]$ as a $C$-set; 
	\item there is an integer $k$ such that for all $a\in A$, $|\{C_0\in \mathfrak C: \exists c\in C_0(a<c)\}|=k$. 
\end{enumerate}
\end{definition}

The last two conditions entail that the $C$-relation cannot distinguish two chains from $\mathfrak C$. 

\begin{remark}\label{finchain}
In the previous proof, the only possible case where we have an infinite chain in the range of $f$ is case 3, i.e., the case where $\phi(M,\alpha)$ is a 1-level set. Furthermore, in order to get the decomposition of the previous theorem, the only new parameters we need to add arise from selecting individual chains from a finite definable set of chains. Notice that using the predicate $C$, one can always definably without parameters partition a finite union of chains into finitely many finite definable sets of incomparable chains. Thus the theorem could also be phrased as: $\dom(f)$ can be decomposed into cells $B_1,\ldots,B_m$ (without adding parameters) such that for each $1\leq i\leq n$, either $f(B_i)$ is an antichain or $f(B_i)$ is the union of a finite set of incomparable chains.  
\end{remark}

The previous proof shows a little bit more:

\begin{corollary}\label{order}
Let $I\subseteq T$ be a chain and $f:M\rightarrow I$ be a partial definable function. Then there are definable chains $B_1,\ldots,B_n$, definable functions $\hat{f_i}:B_i\rightarrow I$ and a partition of $\dom(f)$ into sets $D_1,\ldots,D_{n+1}$ such that $f(D_{n+1})$ is finite and for all $1\leq i\leq n$ if $\alpha\in D_i$ then $f(\alpha)=\hat{f_i}(\inf(\alpha,B_i))$.
\end{corollary}

\pf First, by \ref{monotonicity}, removing finitely many points from $\dom(f)$ we can assume that $f$ is locally constant (if $\dom(f)$ contains a cone where $f$ is a $C$-isomorphism its image would be an antichain, which contradicts that $I$ is totally ordered). Let $g:M\rightarrow T$ be the definable function sending $\alpha\in \dom(f)$ to the basis of the maximal cone containing $\alpha$ on which $f$ is constant. It follows from the previous proof that there is a finite subset $F\subseteq f(M)$ such that for $D:=f^{-1}(F)$ and all $\alpha\in\dom(f)\setminus D$ the set $\Lambda_{g(\alpha)}\cap f^{-1}(f(\alpha))$ is a 1-level set and that $B:=\{g(\alpha)\in T:\alpha \in \dom(f)\setminus D\}$ is a finite union of definable chains $B_1,\ldots,B_n$ (which we may assume without loss of generality disjoint). Set $D_{n+1}$ as $D$ and for $1\leq i\leq n$ define $D_i:= g^{-1}(B_i)$ and the functions $\hat{f_i}: B_i\rightarrow I$ by $\hat{f_i}(b)=a$ if and only if $f(\alpha)= b$ for some (all) $\alpha\in \dom(f)\setminus D_i$ such that $g(\alpha)=b$. These functions are well-defined and satisfy all the required conditions. $\square$

\

This has a direct application to algebraically closed valued fields. Recall that a valued field $K$ is \emph{definably complete} if any definable pseudo-Cauchy sequence has a limit in $K$. 

\begin{corollary}
Let $M$ be a definably complete $C$-minimal expansion of an algebraically closed valued field $K$. We denote by $vK$ the valuation group. Let $f:K\rightarrow vK$ be a partial definable function. Then there are a finite set $F\subseteq f(K)$ and $\beta_1,\ldots,\beta_n$ such that $\dom(f)\setminus f^{-1}(F)$ can be decomposed into cells $D_1\ldots,D_n$ satisfying that for all $1\leq i\leq n$ and all $a\in f(D_i)$ there is $a_i\in vK$ such that
	$$\{x\in D_i: f(x)=a\}=\{x\in D_i: v(x-\beta_i)=a_i\}$$
\end{corollary}

\pf By \ref{monotonicity}, $f$ is locally constant on a cofinite subset of $\dom(f)$, so throwing away finitely many points, we may assume that $f$ is locally constant in $\dom(f)$. For $\alpha\in \dom(f)$, let $g(\alpha)$ be the basis of the maximal cone containing $\alpha$ on which $f$ is constant. Then, by Corollary \ref{order}, there is $F_0$ a finite subset of $f(M)$ such that for all $\alpha\in\dom(f)$, if $f(\alpha)\notin F_0$ then the sets $\Lambda_{g(\alpha)}\cap f^{-1}(f(\alpha))$ are 1-level sets and $B:=\{g(\alpha)\in T:f(\alpha)\notin F_0\}$ is a finite union of chains $B_1,\ldots,B_n$. Let $\beta_1,\ldots,\beta_n$ be leaves above $B_1,\ldots,B_n$ (which exist by the definable completeness assumption). For each $1\leq i\leq n$ we $\beta_i$-definably identify $Br(\beta_i)$ with $vK$ (hence $B_i$ is identified with a subset of $vK$). Let $h_i:vK\rightarrow vK$ be the definable partial function defined on $B_i$ by $h_i(a):=f(\alpha)$ for some (any) $\alpha$ such that $a=g(\alpha)$. By o-minimality, we can partition $\dom(h_i)$ into cells $E_{i1},\ldots,E_{ik_i}$ (o-minimal cells) such that $h_i$ is either constant or monotone on each $E_{ij}$ for $1\leq j\leq k_i$. Let $F_i\subseteq h_i(vK)$ be the union of the sets $h_i(E_{ij})$ for $j$ such that $h_i$ is constant on $E_{ij}$. The set $F=F_0\cup \bigcup_{i=1}^n F_i$ is finite. Let $D_{ij}:=\{\alpha\in \dom(f): g(\alpha)\in E_{ij}\}$ for all $1\leq i\leq n$ and $1\leq j\leq k_i$ such that $h_i$ is not constant (hence is monotone) in $E_{ij}$. By redecomposing into cells, we may assume that $D_{ij}$ is a cell. For $a\in f(D_{ij})\setminus F$ we have that

\begin{center}
	\begin{tabular}{lll}
	$\{x\in D_{ij}: f(x)=a\}$ 	& $=$ & $\{x\in D_{ij}: h_i(g(x))=a\}$ \\
																	& $=$ & $\{x\in D_{ij}: v(x-\beta_i)=h_i^{-1}(a)\}$.\\
\end{tabular}
\end{center}

Setting $a_i$ as $h_i^{-1}(a)$ we obtain what was needed. $\square$

\

As said in the introduction, Corollary \ref{order} shows that the domain of any partial definable function from $M$ to a branch $B$ of $T$ can be decomposed into finitely many cells $D_1,\ldots,D_n$ where the function can be factored (modulo some finite set in the range) as a function from another branch $B'$ to $B$ with the function sending $\alpha\in M$ to its meet with $B'$. We show the analogous result for functions having their range in the set of cones at a point $a\in T$, which we denote by $\mC[a]$. We show first that such functions are also locally constant on a cofinite subset of their domain. 

\begin{lemma}\label{cofini} Let $a\in T$ and $f:M\rightarrow \mC[a]$ be a partial definable function (i.e. for all $\alpha\in \dom(f)$, $f(\alpha)$ is a cone at $a$). Then $f$ is a locally constant function on a cofinite subset of $\dom(f)$. 
\end{lemma} 

\pf Suppose for a contraction that there is a cone $D\subseteq \dom(f)$ such that for all $\alpha\in D$ there is no cone containing $\alpha$ on which $f$ is constant. Therefore, by compactness there is a positive integer $k$ such that for all $\alpha\in D$ we must have that $|f^{-1}(f(\alpha))|<k$. This implies that for all cones $D'\subseteq D$,  $f(D')$ is infinite. But then, taking pairwise disjoint cones $D_1,\ldots,D_k$ in $D$, there must be some $1\leq i\leq k$ such that $f(D_i)$ is an infinite and coinfinite set in $\mC[a]$, which contradicts strong minimality. $\square$

\begin{proposition}\label{setcones}
Let $a\in T$ and let $f:M\rightarrow \mC[a]$ be a partial definable function. Then there are a positive integer $s$, an antichain $\{b_1,\ldots, b_s\}\subseteq T$, a definable set $D\subseteq \dom(f)$ such that $f(D)$ is finite (i.e. $f(D)$ contains finitely many cones at $a$) and definable functions $\hat{f}_{b_i}:\mC[b_i]\rightarrow \mC[a]$ such that $\dom(f)\subseteq \bigcup\{\Lambda_{b_i}:1\leq i\leq n\}$ and if $\alpha\in \Lambda_{b_i}\setminus D$, then $f(\alpha)=\hat{f}_{b_i}(\Gamma_{b_i}(\alpha))$.
\end{proposition}

\pf By Lemma \ref{cofini}, we can assume that $f$ is locally constant. Let $g:M\rightarrow T$ be the definable function sending $\alpha\in \dom(f)$ to the basis of the maximal cone containing $\alpha$ on which $f$ is constant (notice $\dom(f)=\dom(g)$). By Proposition \ref{localcst}, we can decompose the domain of $g$ into finitely many cells $E_1,\ldots, E_n$ such that $g(E_i)$ is either a chain or an antichain for $1\leq i\leq n$. Having the result for each $f\upharpoonright E_i$ for $1\leq i\leq n$ implies the result for $f$, thus we can suppose that $\dom(f)=E_i$ for some $1\leq i\leq n$. 

\

Suppose first that $g(E_i)$ is finite. For $c\in g(E_i)$ let $W_c:=g^{-1}(c)$ and let $I_0:=\{c\in g(E_i): f(W_c) \text{ contains finitely many cones at } a\}$ (which is definable by $C$-minimality). If $I_0=g(E_i)$, then we can set $D:=\dom(f)$ and we are done. Otherwise, let $b_1,\ldots,b_s$ be a list of the elements in $g(E_i)\setminus I_0$. Then the functions $\hat{f}_{b_i}:\mC[b_i]\rightarrow \mC[a]$ defined by $\hat{f}_{b_i}(\Gamma_{b_i}(\alpha)):=f(\alpha)$ for some $\alpha\in \dom(f)\cap \Lambda_{b_i}$ are well-defined. 

\

So suppose from now on that $g(E_i)$ is infinite. We split in two cases depending whether $g(E_i)$ is a chain or an antichain:

\vspace{0.5cm}
\textbf{(1)}	Suppose $g(E_i)$ is an infinite chain $I$. For $\alpha\in E_i$ let $$W(\alpha):=\Lambda_{g(\alpha)}(I)\cap f^{-1}(f(\alpha)).$$ 
Here the notation $\Lambda_{g(\alpha)}(I)$ has the natural meaning, that is, $\Lambda_{g(\alpha)}(I)=\Lambda_{g(\alpha)}(\delta)$ for some (any) $\delta$ such that $I\subseteq Br(\delta)$. By Lemma \ref{fait1}, throwing away finitely many points from $I$ and further decomposing, we may assume that for any $\alpha\in E_i$, $W(\alpha)=\Lambda_{g(\alpha)}(I)$. Then, we get a function $h:I\rightarrow \mC[a]$ such that $h(g(\alpha))=f(\alpha)$. But by strong minimality $h(I)$ must have a finite number of values and we are done (we can include them in $D$).

\vspace{0.5cm}
\textbf{(2)} Suppose that $g(E_i)$ is an infinite antichain $A$. We show this case does not occur. Let $W(\alpha):=\Lambda_{g(\alpha)}\cap f^{-1}(f(\alpha))$. By compactness, Lemma \ref{antichain2} and strong minimality, for all but finitely many $g(\alpha)\in A$, the set $W(\alpha)$ is a 0-level set, so by further decomposing we may assume this is true for all $\alpha\in E_i$. For $c\in A$ let $V_c:=\{b\in A: f(\Lambda_c)=f(\Lambda_b)\}$. By compactness there is an integer $s$ such that $|V_c|\leq s$, for if not, some set $V_c$ would be infinite and contain a cone, which contradicts the definition of $g$. Notice that this implies that $f(\bigcup_{c\in A}\Lambda_c)$ is an infinite union of cones at $a$. Let $A':=cl_{\inf}(A)\setminus A$. By compactness we have two cases 

\begin{enumerate}
	\item[(a)] there is an infinite definable chain $B\subseteq A'$ (possibly adding new parameters);
	\item[(b)] there is a positive integer $k$ such that any totally ordered subset $H$ of $A'$ has cardinality less than $k$.
\end{enumerate}

\

Suppose (a) and for $b\in B$ let $W_b:=\Lambda_{b}(B)\cap \dom(f)$. First consider the subset of $B'\subseteq B$ defined by $B':=\{b\in B: f(W_b) \text{ contains finitely many cones at } a\}$ (which is definable since $bn(a)\geq \aleph_0$). If $B'$ is infinite, the function sending $b\in B$ to the finite set of cones $f(W_b)$ contradicts the stability of $\mC[a]$ as since all $V_c$ are finite, the function must have infinite range. Thus $B\setminus B'$ is infinite. Furthermore, the set $B'':=\{b\in B\setminus B': f(W_b) \text{ contains all cones at } a\}$ is also finite, for if it is infinite there will be some $c\in A$ such that $V_c$ is infinite which is impossible. Thus $B_0:=B\setminus (B'\cup B'')$ is infinite. For each $b\in B_0$ the set $h(b):=\Lambda_a\setminus f(W_b)$ is a finite union of cones at $a$. By strong minimality, there are a positive integer $r$ and an infinite interval $J\subseteq B_0$ such that $h(b)$ is a union of $r$ cones at $a$ and $f(\bigcup_{c\in J}\Lambda_c)$ is an infinite union of cones at $a$. But then the function $h: J\rightarrow \rsub{\mC[a]}{r}$ contradicts the stability of $\mC[a]$. 

\

Thus we may suppose that (b) holds. This implies in particular that every $c\in A$ has a predecessor in $cl_{\inf}(A)$, which contradicts Lemma \ref{denseanti}. $\square$
 
\begin{corollary}
Let $M$ be a $C$-minimal expansion of an algebraically closed valued field $K$. We denote by $K/v$ the residue field. Let $f:K\rightarrow K/v$ be a partial definable function. Then, there are a finite set $F\subseteq f(K)$, a positive integer $s$, elements $\alpha_1,\ldots,\alpha_s,\beta_1,\ldots,\beta_s\in K$ and definable functions $h_{i}:K/v\rightarrow K/v$ for $1\leq i\leq s$ such that for $X:=\dom(f)\setminus f^{-1}(F)$ and $t_i:=v(\alpha_i-\beta_i)$ (for $1\leq i\leq s$)
 $$X\subseteq \bigcup_{1\leq i\leq s} \Lambda_{t_i}(\alpha_i),$$
(identifying $t_i$ with $\inf(\alpha_i,\beta_i)$) and for $1\leq i\leq s$ and $x\in \Lambda_{t_i}(\alpha_i)\cap X$,
$$f(x)=h_{i}\left(\left(\frac{x-\alpha_i}{\alpha_i-\beta_i}\right)/v\right).$$

\end{corollary}

\pf In an algebraically closed valued field each set of cones $\mC[a]$ for $a\in T$ is a copy of the residue field. Therefore we can see $f$ as a function from $K$ to $\mC[\inf(0,1)]$. By \ref{setcones}, there are a positive integer $s$, an antichain $\{t_1,\ldots, t_s\}\subseteq T$, a definable set $D\subseteq \dom(f)$ such that $f(D)$ is finite and definable functions $\hat{f}_{i}:\mC[t_i]\rightarrow \mC[\inf(0,1)]$ such that $\dom(f)\subseteq \bigcup_{1\leq i\leq s} \Lambda_{t_i}$ and if $x\in \Lambda_{t_i}\setminus D$, then $f(x)=\hat{f}_{t_i}(\Gamma_{t_i}(x))$. Set $F$ as $f(D)$, so $X=D$. For $1\leq i\leq s$, we choose representatives $\alpha_i,\beta_i\in K$ such that $\inf(\alpha_i,\beta_i)=t_i$. Without loss of generality, we can assume that for all $1\leq i\leq s$ the open ball centered at $\alpha_i$ of radius $t_i$ is contained in $D$ (hence not in $X$). For $1\leq i\leq s$, let $g_i:\mC[\inf(0,1)]\rightarrow \mC[t_i]$ be the function $g_i(\Gamma_{\inf(0,1)}(x))=\Gamma_{t_i}(x(\alpha_i-\beta_i)+\alpha_i)$, which is well-defined independently of the chosen representative. Then the functions $h_i:\hat{f}_i \circ g_i$ satisfy all the requirements. $\square$

\section{Definable subsets of $T$}

As in the previous section we work in a dense $C$-minimal structure such that for all $a\in T$, $bn(a)\geq\aleph_0$ and for all branches $B\subseteq T$, the order of $B$ is dense. 

\begin{definition}
Let $A\subseteq T$ be a definable antichain and $f:A\rightarrow T$ be a definable function. We say that $f$ is a $t$-function if $f(A)$ is either an antichain or the union of a finite set of incomparable chains.
\end{definition}

\begin{definition}\label{tfun}
A subset $X\subseteq T$ is a 1-cell if there are a definable antichain $A$ and definable $t$-functions $f,g:A\rightarrow T$ such that for all $a\in A$, $g(a)<f(a)$, and $X$ is either $X:=\{t\in T: f(a)=t, a\in A\}$ or $X:=\{t\in T: g(a)<t<g(a), a\in A\}$. 	
\end{definition}

\begin{proposition}\label{tdecomp}
Let $A\subseteq T$ be a definable antichain and $f:A\rightarrow T$ be a definable locally constant function. Then $\dom(f)$ can be definably decomposed (without adding parameters except those used to define $f$) into sets $D_1,\ldots,D_n$ such that for each $1\leq i\leq n$ the function $f\upharpoonright D_i$ is a $t$-function.
\end{proposition}

\pf We may assume that $A$ is infinite. In view of Lemma \ref{denseanti}, we can see $M[A]$ as a $C$-minimal structure. The proposition follows then by the reformulation of Proposition \ref{localcst} in remark \ref{finchain} and the definition of $t$-function. $\square$

\begin{lemma}\label{selflocal}
Let $A$ be a definable antichain and $f:A\rightarrow T$ be a definable function such that for all $a\in A$, $f(a)<a$. Then $f$ is locally constant on a cofinite subset of $A$.
\end{lemma}

\pf As before, we may assume that $A$ is infinite and we can see $M[A]$ as a $C$-minimal structure. By \ref{monotonicity} we can decompose $\dom(f)$ into definable sets $F\cup I\cup K$ where $F$ is finite, both $I$ and $K$ are open, $f$ is a local $C$-isomorphism on $I$ and locally constant on $K$. We show that $I$ is empty. If not, by Lemma \ref{denseanti}, let $a\in I$ and $D\subseteq A$ be a cone containing $a$ on which $f$ is a $C$-isomorphism. Since $f(a)<a$, by density (this is again Lemma \ref{denseanti}) we may suppose that $f(a)<\min(D)$. By definition, $f(D)$ is an antichain, but for $b\in D$ such that $a\neq b$, we have that $f(a)$ is comparable to $f(b)$, a contradiction. $\square$

\

From now on we let $L^*$ be an expansion by definitions of $L$ with 2 sorts, one for the elements of the $C$-set and one for the elements in $T$. 

\begin{proposition}\label{decompT}
Let $\phi(t,y)$ be a $L^*$-formula with $t$ a variable of sort $T$ and $y$ is a tuple of variables of sort $M$. Then there are $n<\omega$ and $L^*$-formulas $\psi_1(t,y),\ldots,\psi_n(t,y)$ such that for all $\alpha\in M^{|y|}$, $\psi_i(T, \alpha)$ is a 1-cell for all $1\leq i\leq n$ and $\phi(T,\alpha)=\bigcup_{i=1}^n \psi_i(T,\alpha)$. Moreover one can assume $\psi_i(T,\alpha)\cap\psi_j(T,\alpha)=\emptyset$ for $1\leq i<j\leq n$.
\end{proposition}

\pf
Fix $\alpha\in M^{|y|}$ and let $X:=\phi(T,\alpha)$. Let $\theta(a):=\exists t \in X (t< a)$ and $\xi(a):=\neg\exists t\in X(a<t)$. Consider the following set
$$A:=\{a\in T\cup M: \theta(a)\wedge \xi(a)\wedge \forall y\in T\cup M((\theta(y)\wedge\xi(y)\wedge a\nparallel y)\rightarrow a\leq y)\}.$$

The set $A$ corresponds to the set of elements in $T$ or $M$ which are bigger than all elements in $X$ to which they can be compared with and are minimal for that property. It is not difficult to see that $A$ is a definable antichain. For $a\in A$, let $Y_a:=\{t\in X: t\leq a\}$. By compactness and o-minimality, there are positive integers $k$ and $l$ such that for all $a\in A$, $Y_a$ (viewed as a subset of $\{t\in T:t\leq a\}$) is a union of at most $k$ infinite open intervals and $l$ isolated points. We partition $A$ into finitely many parts such that $a,b\in A$ belong to the same equivalence class if and only if $Y_a$ and $Y_b$ have the same number of open intervals and the same number of isolated points. Then refine the partition such that the order of these intervals and points is the same (for example, if each $Y_a$ has two intervals and an isolated point, then there are three ways of ordering them, either the point is in between the two intervals, above them or below them). We may suppose then that $A$ is one of these equivalence classes, i.e., a definable antichain such that for all $a\in A$, $Y_a$ has exactly $k$ open intervals and $l$ isolated points which are ordered in the same way for each $a\in A$. Take now $X$ as $X_A:=\{t\in X: t\leq a, \text{ for some } a\in A\}$. For $1\leq i\leq k$, let then $f_i:A\rightarrow T$ be the definable functions sending $a$ to the right endpoint of the $i^\text{th}$ interval in $Y_a$, $g_i:A\rightarrow T$ be the definable functions sending $a$ to the left endpoint of the $i^\text{th}$ interval in $Y_a$ and for $1\leq i\leq l$ let $h_i:A\rightarrow T$ be the definable functions sending $a$ to the $i^\text{th}$ isolated point in $Y_a$. We decompose $X$ into sets $X_1,\ldots,X_k$ and $W_1,\ldots,W_l$ where

\begin{center}
	\begin{tabular}{l}
$X_i:=\{t\in X: g_i(a)<t<f_i(a), \text{ for any } a\in A\}$ and \\
$W_i:=\{t\in X: h_i(a)=t, \text{ for any } a\in A\}$.\\
\end{tabular}
\end{center}

By Lemma \ref{selflocal}, all functions $f_i,g_i$ and $h_i$ are locally constant on a cofinite subset of $A$. Thus removing finitely many points we can suppose they are locally constant on $A$ (notice removing finitely many points can always be done for our purposes). By applying Proposition \ref{tdecomp} repeatedly, the functions $f_i,g_i$ and $h_i$ decompose $A$ into finitely many pieces such that on each piece each function $f_i,g_i$ and $h_i$ is a $t$-function. This gives us a union of 1-cells that is equal to $X$. To get a disjoint union, one has to take the minimal refinement of such a union. Notice that through the proof all definable sets used are uniformly definable in $\alpha$. $\square$

\bibliographystyle{jflnat}
\bibliography{models}

\begin{thebibliography}{5}
\expandafter\ifx\csname natexlab\endcsname\relax\def\natexlab#1{#1}\fi
\def\docolon{:}
\def\eatcomma#1{}
\def\onlyone#1{\gdef\oneletter{#1}}
\def\sphref#1#2{{\let\#=\docolon\xdef\one{#1}}\href{\one}{#2}}
\def\zhref#1,#2{{\let\#=\docolon\xdef\one{#1}}\href{\one}{#2}}
\expandafter\ifx\csname url\endcsname\relax
  \def\url#1{{\tt #1}}\fi
\newcommand{\enquote}[2]{``#1,''}

\bibitem[{Cubides-Kovacsics}(2014)]{cubides:12}
{Cubides-Kovacsics}, P.,
\newblock
\newblock \enquote{An introduction to {$C$}-minimal structures and their cell
  decomposition theorem}, in {\em Valuation Theory in Interaction}, edited by
  \reverseeditorsnames Antonio Campillo and Franz-Viktor Kuhlmann and Bernard
  Teissier, Series of Congress Reports of the EMS, pp.~167--207. European
  Mathematical Society, 2014.

\bibitem[Delon(2011)]{delon:11}
Delon, F.,
\newblock \enquote{C-minimal structures without the density
  assumption}\eatcomma,
\newblock  in {\em Motivic Integration and its Interactions with Model Theory
  and Non-Archimedean Geometry}, edited by \reverseeditorsnames Raf Cluckers
  and Johannes Nicaise and Julien Sebag,
\newblock Cambridge University Press, Berlin, 2011.

\bibitem[Haskell et~al.(2008)Haskell, Hrushovski, and
  Macpherson]{haskellETAL:08}
Haskell, D., E.~Hrushovski \unskip, and  D.~Macpherson,
\newblock {\em Stable domination and independence in algebraically closed
  valued fields}, volume~Lecture Notes in Logic,
\newblock Cambridge Unversity Press, 2008.\eatcomma.

\bibitem[Haskell and Macpherson(1994)]{macphersonETAL:94}
Haskell, D. \unskip, and  D.~Macpherson,
\newblock \enquote{Cell decompositions of {C}-minimal structures},
\newblock {\em Annals of Pure and Applied Logic}, vol.~66 (1994),
  pp.~113--162.\eatcomma.

\bibitem[Macpherson and Steinhorn(1996)]{macphersonETAL:96}
Macpherson, D. \unskip, and  C.~Steinhorn,
\newblock \enquote{On variants of o-minimality},
\newblock {\em Ann. Pure Appl. Logic}, vol.~79 (1996), pp.~165--209.\eatcomma.

\end{thebibliography}

\end{document}